\documentclass{elsart3-1VARIANT} 
\usepackage{amssymb}
\usepackage[english,francais]{babel} 
\usepackage[mathcal]{euscript}
\usepackage{amstext}
\usepackage{amssymb,mathrsfs} 

\newtheorem{theorem}{Theorem}[section]

\newtheorem{theoreme}{Th\'eor\`eme}[section]

\newtheorem{definition}[theoreme]{Definition}   

\newcommand \omegab 	{\overline \omega}
\newcommand \iint 	{\int\hskip-.5cm\int}
\newcommand \la 		\langle
\newcommand \ra 		\rangle
\newcommand \be   	{\begin{equation}}
\newcommand \ee   	{\end{equation}}
\newcommand \RR      {\mathbb{R}}
\newcommand \eps     \epsilon
\newcommand \ubar   	{{\overline u}} 
\newcommand \dive 	{\mbox{div}} 
\newcommand \del		\partial  
\newcommand \Ecal {\mathcal E} 

\newcommand \lad 		{\prec}
\newcommand \rad 		{\succ} 
\journal{Comptes Rendus de l'Acad\'emie des Sciences}
\begin{document} 
 
\begin{frontmatter}
 
\selectlanguage{english}
\title{Hyperbolic conservation laws on manifolds
\\
 with limited regularity
 }

\selectlanguage{english}
\author[authorlabel]{Philippe G. LeFloch},
%\ead{LeFloch@ann.jussieu.fr}
\author[authorlabel]{Baver Okutmustur}
%\ead{Okutmustur@ann.jussieu.fr}

\address[authorlabel]{Laboratoire Jacques-Louis Lions, Centre National de la Recherche Scientifique
\\
Universit\'e Pierre et Marie Curie (Paris 6)
\\4 Place Jussieu,  75252 Paris, France.
\\
{\sl E-mail: \tt LeFloch@ann.jussieu.fr, Okutmustur@ann.jussieu.fr} 
}

\medskip
\begin{center} 
\end{center}

\begin{abstract}
\selectlanguage{english}
We introduce a formulation of the initial and boundary value problem for nonlinear hyperbolic conservation laws posed  
on a differential manifold endowed with a volume form, possibly with a boundary; 
in particular, this includes the important case of Lorentzian manifolds.
Only limited regularity is assumed on the geometry of the manifold. 
For this problem, we establish the existence and uniqueness of an $L^1$ semi-group 
of weak solutions satisfying suitable entropy and boundary conditions.

\vskip 0.5\baselineskip

\selectlanguage{francais}
\noindent{\bf R\'esum\'e} \vskip 0.5\baselineskip \noindent
{\bf Lois de conservation hyperboliques sur les vari\'et\'es \`a faible r\'egularit\'e. }
Nous proposons une formulation du probl\`eme de Cauchy avec conditions aux limites 
pour les lois de conservation hyperboliques
nonlin\'eaires pos\'ees sur une vari\'et\'e diff\'erentiable munie d'une forme volume, 
avec ou sans bord; notre \'etude couvre, en particulier, le cas important des vari\'et\'e Lorentzienne. 
Nous supposons une r\'egularit\'e limit\'ee sur la g\'eometrie de la vari\'et\'e. 
Pour ce probl\`eme nous d\'emontrons l'existence et l'unicit\'e 
d'un semi-groupe $L^1$ de solutions faibles satisfaisant \`a des conditions d'entropie et 
 \`a des conditions aux limites convenablement d\'efinies. 

\

Reference of this paper: {\tt C.R. Math. Acad. Sc. Paris 346 (2008), 539--543.}

\end{abstract}
\end{frontmatter}

\selectlanguage{francais}
\section*{Version fran\c{c}aise abr\'eg\'ee}

Dans cette Note, nous nous int\'eressons aux lois de conservation hyperboliques nonlin\'eaires pos\'ees 
sur une vari\'et\'e. Nous proposons une formulation math\'ematique du probl\`eme de Cauchy 
lorsque la vari\'et\'e est seulement munie d'une forme volume. 
En particulier, nous pouvons poser ce probl\`eme sur une vari\'et\'e Lorentzienne. 

Soit $M$ une vari\'et\'e compacte, r\'eguli\`ere, \`a $n$ dimensions, 
munie d'une forme volume $\omega$ de class $L^\infty$. Sur cette vari\'et\'e, consid\'erons la loi de conservation 
hyperbolique nonlin\'eaire 
\be
\label{CL.1f}
  \partial_t u + \dive_\omega \big( f(u) \big) = 0,
  \qquad u: \RR_+\times M \to \RR, 
\ee
o\`u le flux $f=(f^j(\ubar,x))$ est un champ de vecteur r\'egulier d\'efini sur $M$ 
et d\'ependant d'un param\`etre r\'eel $\ubar$. L'op\'erateur de divergence spatial est ici
 d\'efini au sens faible des distributions. (Voir (\ref{star}) dans la Version en anglais.) 

Nous supposons que le flux $f$ est au plus lin\'eaire, c'est-\`a-dire que  
pour tout champ $\alpha$ de $1$-formes de classe $L^\infty$ d\'efinies sur $M$ il existe une constante 
$C_\alpha>0$ telle que 
\be
\label{growth0}
\sup_M \left| \la \alpha, f (\ubar) \ra \right| 
\leq C_\alpha \, (1 + |\ubar|), \qquad \ubar \in \RR.   
\ee
Nous recherchons alors des solutions verifiant une condition initiale en $t=0$:
\be
\label{CL.2f}
u(0,x) = u_0(x), \qquad x\in  M,
\ee
o\`u $u_0 \in L^1_\omega(M)$ est donn\'ee. Nous disons que $f$ est compatible avec la g\'eom\'etrie de $M$
si 
\be
\label{CL.3f}
\big( \dive_\omega f\big) (\ubar) = 0, \qquad \ubar \in \RR.
\ee
Nous d\'efinissons alors une notion de couple d'entropie convexe (U,F) adapt\'ee \`a la vari\'et\'e $M$. 
(Voir (\ref{starstar}) dans la Version en anglais.) 
\

\begin{definition} 
Sous les hypoth\`eses ci-dessus, une fonction $u \in L^\infty (\RR_+, L^1_\omega( M))$ est appel\'ee
solution entropique du probl\`eme (\ref{CL.1f})-(\ref{CL.2f})
si pour tout couple d'entropie convexe $(U,F)$ et toute fonction-test 
$\phi : [0,\infty) \times M \to \RR_+$,
\begin{eqnarray*}
&&\iint_{\RR_+ \times M} \Big( U(u) \, \del_t \phi +  \la d\phi, F(u)\ra \Big) \, \omega dt
   + \iint_{\RR_+ \times M} \Big( \big( \dive_\omega F\big)(u) - \del_u U(u) \big(\dive_\omega f\big)(u) \Big) \, \phi \, \omega dt
\\
&& + \int_M U(u_0) \phi(0) \, \omega \geq 0.
\end{eqnarray*}
\end{definition}

\begin{theorem}
\label{CL-0}
Soit $M$ une vari\'et\'e compacte sans bord et $\omega \in L^\infty(M)$ 
une forme volume uniform\'ement d\'efinie positive 
sur $M$. Soit $f$ un flux compatible avec la g\'eom\'etrie de $M$ et 
satisfaisant \`a la condition de croissance (\ref{growth0}). 
Alors, le probl\`eme (\ref{CL.1f})--(\ref{CL.2f}) admet un semi-groupe contractif de solutions 
entropiques  $u_0 \in L_\omega^1(M) \mapsto u(t) := S_t u_0 \in L_\omega^1(M)$
$$
\| S_t u_0 - S_t v_0 \|_{L^1_\omega(M)} 
\leq
\|  u_0 -  v_0 \|_{L_\omega^1(M)},
\qquad t \geq 0, \quad u_0, v_0 \in L_\omega^1(M).
$$
\end{theorem}

\

Nous \'etendons aussi ce r\'esultat au cas des lois de conservation d\'efinies sur 
un espace-temps muni d'une m\'etrique Lorentzienne de classe $L^\infty$, 
et incluons des conditions aux limites lorsque $M$ est une 
vari\'et\'e \`a bord. Lorsque la condition de compatibilit\'e g\'eom\'etrique est relax\'ee, 
la distance $L^1$ entre deux solutions entropiques peut cro\^{\i}tre en temps.
Pour plus de d\'etails nous renvoyons le lecteur \`a la version en anglais et \`a \cite{LO}.  

%====================================================================================================================

\selectlanguage{english}
\section*{Version in English}
\section{Introduction}
\label{IN}

In this Note we consider nonlinear hyperbolic conservation laws posed on a differentiable manifold
and we introduce a formulation of the initial and boundary value problem. 
The manifold under consideration 
may have a boundary and is endowed with a Lorentzian metric or, more generally, a volume form. 
We introduce a notion of weak solution in the sense of distributions and 
formulate suitable entropy and boundary conditions. 
We then establish the existence and uniqueness of entropy solutions when the initial data are integrable 
functions and the manifold geometry has limited regularity. 
This problem is motivated by similar questions arising in the evolution of compressible fluids 
on manifolds (e.g. fluid flows on the sphere or on an Einstein spacetime of general relativity).
Scalar conservation laws provide us with a simplified model for understanding certain important 
features arising in this context, especially the lack of regularity of the solutions and the geometry. 

Our results generalize the well-posedness and convergence theory established in \cite{ABL,BL,BFL,PLF,Panov}  
for the case of Riemannian manifolds and, earlier on, in Kruzkov's classical work \cite{Kruzkov} in the Euclidian case. 
Our presentation will proceed by investigating three 
different settings that are (related but) of independent interest in the applications: 
$n$-dimensional manifolds endowed with a volume form, Lorentzian manifolds, 
and $(n+1)$-dimensional manifolds endowed with a volume form. 
Our proofs use a sequence of approximate solutions constructed by the finite volume method.

%===============================================================================================================

\section{Conservation laws on a manifold endowed with an $L^\infty$ volume form}

We consider first the case that $M$ is a compact, smooth, $n$-dimensional manifold
endowed with a volume form $\omega \in L^\infty(M)$. So, $\omega$ is a bounded Lebesgue measurable, uniformly
positive $n$-form field which, in local coordinates $x=(x^j)_{1 \leq j \leq n}$  defined in some open set $B$, reads 
$$
\omega = \omegab \, dx^1 \ldots dx^n, \qquad \omegab \in L^\infty(B), \quad \omega \geq \omega_*
$$
for some $\omega_*>0$ (depending on the choice of local coordinates). 
On this manifold, we consider the {\sl hyperbolic conservation law} 
\be
\label{CL.1}
  \partial_t u + \dive_\omega \big( f(u) \big) = 0,
  \qquad u: \RR_+\times M \to \RR, 
\ee
in which the flux $f=(f^j(\ubar,x))$ is a smooth vector field on $M$ depending upon a real parameter $\ubar$.
In (\ref{CL.1}), the spatial divergence operator is defined in a weak sense from the volume form $\omega$. 
Namely, when $u$ and $\omega$ are smooth and the vector field $f(u)$ is supported in a local chart of coordinates,
one defines 
$$
(\dive_\omega f(u))(t,x) : = {1 \over \omegab(x)} \, \del_j \Big( \omegab(x) \, f^j(u(t,x),x) \Big),
$$
where implicit summation over repeated indices is used. Hence, for every smooth function $\theta: M \to \RR$ 
compactly supported in a local coordinate chart and for $t\geq 0$, 
\begin{eqnarray}
\label{star}
\int_M \dive_\omega \big( f(u) \big) (t,\cdot) \, \theta \, \omega
&& = - \int_B (\del_j\theta)(x) \, f^j(u(t,x),x)\, \omegab(x) dx^1 \ldots dx^n
\\
\nonumber 
&& = \int_M \la d\theta, f(u(t, \cdot)) \ra \, \omega, 
\end{eqnarray}
where the one-form field $d\theta$ is the differential of the function $\theta$
and $\la \cdot, \cdot \ra$ denotes the bracket between forms and vectors on $M$.  
The above formula allows one to reformulate the conservation law (\ref{CL.1}) in a {\sl weak sense} 
when $\omega \in L^\infty(M)$, 
$u \in L_\omega^1(M)$, and the flux $f$ is at most linear, i.e., for every 
$L^\infty$ field $\alpha$ of $1$-forms on $M$ there exists a constant $C_\alpha>0$ such that 
\be
\label{growth}
\sup_M \left| \la \alpha, f (\ubar) \ra \right| 
\leq C_\alpha \, (1 + |\ubar|), \qquad \ubar \in \RR.   
\ee

We are interested in solutions satisfying the following initial condition at the time $t=0$:
\be
\label{CL.2}
u(0,x) = u_0(x), \qquad x\in  M,
\ee
where the function $u_0 \in L^1_\omega(M)$ is given. Generalizing \cite{BL},  
we say that $f$ is geometry-compatible if 
\be
\label{CL.3}
\big( \dive_\omega f\big) (\ubar) = 0, \qquad \ubar \in \RR.
\ee
By definition, a Lipschitz continuous function $U:\RR \to \RR$ together with a family of 
vector fields $F=F(\ubar)$ depending Lipschitz continuously upon the parameter $\ubar$
is called an {\sl entropy pair} for the conservation law (\ref{CL.1}) if for almost all $\ubar \in \RR$
\be
\label{starstar}
\del_u F(\ubar) = \del_u U(\ubar) \del_u f(\ubar).
\ee
Finally, we introduce the notion of entropy solution. 

\

\begin{definition} 
Under the above assumptions, a function $u \in L^\infty (\RR_+, L^1_\omega( M))$ is called
an {\rm entropy solution} to the initial value problem (\ref{CL.1})-(\ref{CL.2})
if for every entropy pair $(U,F)$ and every smooth function
$\phi \geq 0$ with compact support in $[0,\infty) \times M$,
\begin{eqnarray*}
&&\iint_{\RR_+ \times M} \Big( U(u) \, \del_t \phi +  \la d\phi, F(u)\ra \Big) \, \omega dt
   + \iint_{\RR_+ \times M} \Big( \big( \dive_\omega F\big)(u) - \del_u U(u) \big(\dive_\omega f\big)(u) \Big) \, \phi \, \omega dt
\\
&& + \int_M U(u_0) \phi(0) \, \omega \geq 0.
\end{eqnarray*}
\end{definition}

\begin{theorem}
\label{CL-02}
Let $M$ be a compact, smooth manifold without boundary and $\omega \in L^\infty(M)$ be a uniformly positive  
volume form on $M$. Let $f$ be a geometry-compatible flux on $M$ satisfying the growth condition (\ref{growth}).
Then, the initial value problem (\ref{CL.1})--(\ref{CL.2}), 
admits a contractive semi-group of entropy solutions $u_0 \in L_\omega^1(M) \mapsto u(t) := S_t u_0 \in L_\omega^1(M)$, 
$$
\| S_t u_0 - S_t v_0 \|_{L^1_\omega(M)} 
\leq
\|  u_0 -  v_0 \|_{L_\omega^1(M)},
\qquad t \geq 0, \quad u_0, v_0 \in L_\omega^1(M).
$$
\end{theorem}

\

The above definition and theorem extend a classical result by Kruzkov \cite{Kruzkov}, which was restricted 
to equations posed on the (flat) Euclidian space. 
The proof follows the lines of arguments in \cite{BL} and relies on the observation that
the volume form structure $\omega$ only is required to carry out all of the arguments therein. 
Indeed, the metric can be eliminated and replaced by the corresponding expressions in terms of the volume form.
Note also that boundary conditions can also be included if $M$ has a boundary 
and the condition (\ref{CL.3}) can also be relaxed; see \cite{LO} for further details.

%==================================================================================================================

\section{Conservation laws on a spacetime endowed with a Lorentzian metric}

We now discuss a setting where further geometric structure is provided on the manifold;  
this is motivated by problems arising in the theory of general relativity. For simplicity
in the presentation,  
we assume in this section that the geometry is sufficiently smooth, and we refer the following section for the case 
of a limited regularity. 
We assume that a smooth, $(n+1)$-dimensional manifold $M$ with smooth boundary $\del M$ 
is endowed with a smooth Lorentzian metric $g$ with signature $(-, +, \ldots, +)$. We denote by $dv_M $ the volume form on $(M,g)$. We say, in short, that $(M,g)$ is a {\sl spacetime.}

A flux on the manifold $M$ is still defined as a smooth vector field $f=f(\ubar,x)$ depending on a parameter $\ubar$,
but it is important to observe that time and space, now, are handled together. So, in local coordinates we may
write $x=(x^\alpha)=(t, x^j)$, with $\alpha=0,\ldots, n$ and $j=1, \ldots, n$, 
and we set $f = \big(f^j(\ubar, x)\big). $
We consider the following conservation law posed on $M$: 
\be
\label{CL.4}
\dive_g \big( f(u) \big) = 0, \qquad u: M \to \RR.
\ee

To formulate the initial and boundary value problem associated with (\ref{CL.4}),
we prescribe a measurable and bounded function $u_B: \del M \to \RR$ defined on the boundary $\del M$ of the spacetime, 
and we search for a function $u \in L^\infty(M)$ satisfying (\ref{CL.4}) in the distributional sense,
 together with 
an entropy condition, such that the (weak) trace of $u$ on $\del M$ satisfies the boundary condition 
\be
u \big|_{\del M} \in \Ecal_N ( u_0)
\label{CL.5}
\ee
in a sense specified now. 
By definition, a convex function $U:\RR \to \RR$ and a vector field $F=F(\ubar)$ is called a convex entropy pair
associated with the conservation law (\ref{CL.4}) if, in local coordinates,  
$$
F^j(\ubar, x) = \int_0^\ubar \del_u U(u') \,\del_u f^j(u',x) \, du', \qquad \ubar \in \RR.
$$
Here,  
$N$ denote the field of unit normal $1$-forms along the boundary, and 
for all convex entropy pair $(U,F)$ we set 
$$
\Ecal_N (u_0):= \Big\{ \ubar \, \, \big| \, \,  
E_N(u_B, \ubar) := \la N, F(u_B) \ra + \del_u U(u_B) \la N, f(\ubar) - f(u_B) \ra 
\leq \la N, F(\ubar) \ra \Big \}.
$$
In the boundary condition (\ref{CL.5}) we do not distinguish between space-like parts (initial hypersurface, final hypersurface) 
or time-like parts of the boundary $\del M$. If $S \subset \del M$ is 
a space-like hypersurface then the boundary condition 
reduces to either a vacuous requirement or else to the continuity property $u\big |_S = u_B$.

\

\begin{definition}
\label{CL-12} 
A Young measure $\nu: M \to \text{Prob}(\RR)$ with compact range is called an {\rm entropy measure-valued solution}  
to the problem (\ref{CL.4})--(\ref{CL.5}) if for some scalar field $b\in L^\infty(\del M)$  
the inequalities  
\begin{eqnarray*}
& \int_M 
\lad \nu, \la d\theta, F\ra + \big( \dive_g F - \del_u U \big(\dive_g f\big) \big) \theta \rad \, d v_M 
+ \int_{\del M} \lad \nu, E_N(u_B, b) \rad \, \theta \, dv_{\del M} \geq 0
\end{eqnarray*}
hold for all convex entropy pairs $(U,F)$ and all smooth functions $\theta \geq 0$ compactly supported in $\overline M$.
On the other hand, a function $u \in L^{\infty}(M)$ is called an {\rm entropy solution}
 to the problem (\ref{CL.4})--(\ref{CL.5})
if and only if the associated Young measure $\delta_u$ 
(the Dirac measure at the point $u$)
is an entropy measure-valued solution to the same problem. 
\end{definition}

\

The above definition extends a notion introduced by DiPerna \cite{DiPerna}
for equations posed in the (flat) Euclidian space and by Szepessy \cite{Szepessy} and 
Kondo and LeFloch \cite{KL} for the problem with boundary conditions. 
The formulation and the convergence of finite volume schemes on a Lorentzian manifold without boundary 
was treated in \cite{ALO}. 
In the following theorem we suppose that $M$ is foliated by oriented space-like hypersurfaces, that is,  
$M = \bigcup_{t\geq 0} H_t$.

\

\begin{theorem}
\label{CL-1}
Suppose that $(M,g)$ is a time-oriented, smooth, $(n+1)$-dimensional Lorentzian manifold foliated by 
oriented space-like hypersurfaces $H_t$. 
Let $f=f(\ubar)$ be a future-oriented, time-like, smooth vector field depending on a parameter. 
Then, for each data $u_B \in L^\infty(\del M) \cap L_g^1(\del M)$,
the initial and boundary value problem (\ref{CL.4})--(\ref{CL.5}) admits 
a unique entropy solution $u= S_t u_B \in L^\infty(M)$ such that $u_{H_t} \in L_g^1(H_t)$ for all times $t \geq 0$, 
with moreover for each $T>0$ and $t \in [0,T]$ 
$$
\| S_T u_0 - S_T v_0 \|_{L^1_g(H_t)} 
\leq  
C_T \, \| u_B - v_B\| _{L^1_g((\del M)_T)},
\qquad u_B, v_B \in L_g^1(\del M)\cap L^\infty(\del M), 
$$
where $(\del M)_T$ is the part of the boundary corresponding to 
$\big\{0\leq t \leq T\big\}$ 
and the constant $C_T$ depends on $T$ and the sup norm of the data, only. 
\end{theorem}

\

We can show that, in the integral term $\| u_B - v_B\| _{L^1_g((\del M)_T)}$, one can suppress the part along 
which the vector field $\del_u f(\ubar, \cdot)$ is outgoing for all~$\ubar$. 
 
%==================================================================================================================

\section{Conservation laws on a spacetime endowed with an $L^\infty$ volume form}

Finally, the above results are extended to the situation where $M$ is an 
$(n+1)$-dimensional spacetime endowed with a volume form $\omega$ (rather than a Lorentzian metric)
and, moreover, is solely of class $L^\infty$. Given a parameter-dependent 
vector field $f$ with ``at most linear'' growth, weak solutions in $L^1_{\text{loc}}(M)$ 
to the hyperbolic conservation law 
\be
\label{CL.4b}
\dive_\omega \big( f(u) \big) = 0, \qquad u: M \to \RR
\ee
are defined 
by adapting the definition introduced in Section~2. Observe that 
the (spacetime) volume form need not admit a (strong) trace on $n$-dimensional hypersurfaces in $M$, so that 
the conclusions of Theorems~\ref{CL-02} and \ref{CL-1} hold in a weak sense only. We refer to \cite{LO} for further details. 

%------------------------------------------------------------------------------------------------------------------------

\end{document}